\documentclass[11pt,a4paper,draft]{article}
\usepackage{amsmath,mathrsfs,amssymb,amsthm,amscd}
\usepackage[]{fontenc}
\usepackage{xy}
\usepackage{enumerate}
\xyoption{all} \numberwithin{equation}{section}
\newtheorem{theorem}[equation]{Theorem}
\newtheorem{corollary}[equation]{Corollary}
\newtheorem{proposition}[equation]{Proposition}
\newtheorem{lemma}[equation]{Lemma}

\theoremstyle{definition}
\newtheorem{definition}[equation]{Definition}

\newtheorem{example}[equation]{Example}

\newtheorem{remark}[equation]{Remark}

\title{On the topological index of irregular surfaces}
\author{M.A. Barja, J.C. Naranjo and G.P. Pirola
\footnote{\noindent M.A.Barja: Partially supported by the Proyecto
de Investigaci\'on BFM2003-06001 and ACI-2003 \newline
J.C.Naranjo: Partially supported by the Proyecto de
Investigaci\'on BFM2003-02914 and ACI-2003 \newline G.P. Pirola:
Partially supported by 1) PRIN 2003  ``Spazi di moduli e teoria di
Lie"; 2) Gnsaga; 3) Far 2004 (PV):
 Variet\`{a} algebriche, calcolo algebrico, grafi orientati e topologici}
}

\begin{document}

\maketitle

\begin{abstract}
We study the topological index of some irregular surfaces that we
call generalized Lagrangian. We show that under certain hypotheses
on the base locus of the Lagrangian system the topological index
is non-negative. For the minimal surfaces of general type with
$q=4$ and $p_g=5$ we prove the same statement without any
hypothesis.
\end{abstract}

\section{Introduction.}

We consider complex smooth algebraic varieties $X$
 of dimension $n$ equipped with an element $\omega $ in the kernel of the map
$$
\psi_{2,\Omega ^1_X}=\psi_2:\Lambda ^2 H^0(X, \Omega
^1_X)\longrightarrow H^0(X,\Omega ^2_X).
$$
Let $V\subset H^0(X,\Omega ^1_X)$ be such that $w\in \Lambda
^2\,V$ and $V$ has minimal dimension with this property. Then $V$
has even dimension $2k$ and we say that this is the \textit{rank
of $w$}. Notice that rank $2$ means decomposable.

\vskip 5mm The non-triviality of the kernel of $$ \Lambda
^2H^1(X,\mathbb C)\longrightarrow H^2(X,\mathbb C) $$
 determines
the existence of nilpotent towers in the fundamental group of the
variety (see \cite{Morgan} and also \cite{Am:wb}). Hence, we
obtain topological consequences  of the non-injectivity of $\psi _
2$. This is especially clear when an element of rank $<2n$ appears
in $Ker(\psi _{2})$; then, using the results of Catanese in
\cite{Cat2}, one sees that $X$ is fibred on a variety of lower
dimension. If $n=2$ this is the classical Castelnuovo-de Franchis
Theorem. Our principal interest is to find numerical restrictions
on the topological invariants of the varieties (mainly surfaces)
where $Ker(\psi _2)\ne 0$.

\vskip 5mm
 Consider the image of the evaluation map on $V$:
 $$
 V\otimes \mathcal {O}_X \longrightarrow  \Omega ^1_X.
 $$
 The image defines a torsion free sheaf $\overline{\Omega ^1_X}$.

 \vskip 5mm
 \begin{definition}
We say that an $n$-dimensional variety $X$ is \textit{generalized
Lagrangian } if there exists $w\in Ker(\psi _2)$ of rank $2n$ and
moreover $rank \, \overline{\Omega^1_X}=n.$ In other words, there
exist $\omega _1,\dots ,\omega _{2n} \in H^{0}(X,\Omega ^1_X)$
generating generically $\Omega ^1_X$ and such that $\omega _1
\wedge \omega _2 +\dots +\omega _{2n-1} \wedge \omega _{2n}=0$ in
$H^0(X,\Omega^2_X).$
\end{definition}

Following \cite{Bog:Tsc}, we say that $X$ is \textit{Lagrangian}
if there exists a degree $1$ map
$$
a:X\longrightarrow a(X) \subset A,
$$
where $A$ is an Abelian variety of dimension $2n$, and  a $(2,0)$
form $w$ of maximal rank $2n$ on $A$ such that $a^{\ast}(w)=0$.
This is a particular case of our definition 1.1.

\vskip 5mm  We notice that in this case the vector space $V$
appears as the pull-back of the cotangent space of $A$.

\vskip 5mm We want to give restrictions on the invariants of these
varieties, mainly regarding their topological significance. We
obtain results of two types: first we work with generalized
Lagrangian varieties of any dimension and we prove that under some
hypotheses the second degree part of the Chern character is
non-negative. In fact, this is a consequence of a more general
result on reflexive sheaves (see \ref{thmgen}).  Second we focus
on the $2$-dimensional case and  we assume that $X=S$ is minimal
of general type.
 It is natural to ask what kind of geometry appears if
 a relation of type
 $$  \omega _1\wedge \omega_2+\omega _3\wedge \omega_4=0, $$ holds.
In other words, the study of generalized Lagrangian surfaces
corresponds naturally to a higher rank analogon  of the classical
situation elucidated by  Castelnuovo and de Franchis. In fact, we
are mainly interested in the  non fibred generalized Lagrangian
surfaces.

\vskip 5mm Attached to the vector space $V$ we have a natural
subsystem of the canonical system: the one given by the image of
$\Lambda ^2V$ in $H^0(S,\omega_S)$. Let $F_V$ be the corresponding
base divisor. Assume that $F_V$ is reduced. We say that $F_V$ is
contracted by $V$ if the map
$$
V\otimes  \mathcal O_S\longrightarrow \Omega^1_S \longrightarrow
\Omega ^1_{S_{|F_V}}\longrightarrow \omega
 _{F_V}
$$
vanishes. Then we prove (see \ref{thmsurFzero} and
\ref{thmsurFnc}):

\begin{theorem}\label{thmintro}
Assume that $F_V=0$ or $F_V$ is a reduced connected divisor with
normal crossings
 contracted by $V$. Then the topological index $\tau(S)$
satisfies $\tau (S)\ge 0$.
\end{theorem}

\vskip 5mm We will see that the condition ``$F_V$ contracted by
$V$" is necessary in general (see section 6). It is natural to
state:

\textbf{Conjecture $1$:} Let $S$ be a minimal generalized
Lagrangian surface of general type. Assume that $F_V$ is
contracted by $V$. Then $\tau (S)\ge 0$.

\textbf{Conjecture $2$:} Let $S$ be a Lagrangian surface. Then
$\tau (S)\ge 0$.

\vskip 5mm Observe that, by \ref{thmintro}, Conjecture 2 is true
when $a:X\longrightarrow A$ is an immersion, and more generally
when $a$ does not have a branch divisor.

\vskip 5mm Notice that the existence of non decomposable $w$
implies $q(S)\ge 4$. As a test case to our conjecture, the first
values where non-fibred surfaces
 could appear are $ q(S)=4 $ and $ p_g(S)=5 $.
 In this case the map $\psi_2 $
 has non-trivial kernel for dimensional reasons.
 We see that with these invariants no hypothesis is needed to obtain the inequality $\tau(S)\ge0$.
 In fact, since $K^2_S+c_2(S)=24$ and $\tau(S)=\frac 13 (K_S^2-2c_2(S))$, we can state our result
as follows:

\begin{theorem}\label{main}
Let $S$ be a minimal smooth projective surface of general type
with $p_g(S)=5$ and $q(S)=4$. Then the following inequality holds:
$$
K_S^2\ge 8\,\chi (\mathcal O_S)=16.
$$
If $S$ has an irreducible pencil, then $K_S^2=16$ and we have a
complete classification of these surfaces (see \ref{fibred}).
\end {theorem}

\vskip 5mm Recall the Bogomolov-Miyaoka-Yau inequality $K_S^2\le 9
\chi (\mathcal O_S)=18$, therefore the only possible values are
$16,17$ and $18$.

\vskip 5mm We believe that this is in itself an interesting
classification result on irregular surfaces. The result
\ref{thmintro} is crucial in order to show it. The proof is rather
involved and uses counting quadrics containing the canonical
image,
  an intensive use of Reider's Theorem (\cite{Reider}) and of the Hodge-Index theorem.
 The method exploits heavily the particular geometry of the surfaces with these invariants and
 does not seem suited to approaching the general situation posed in our conjecture.

\vskip 5mm
 The paper is organized as follows: in section 3
 we study the kernel of the map $\psi_{2,\mathcal E}$
 for any reflexive sheaf $\mathcal E$.
 At the end of the section we find the results on generalized Lagrangian varieties of any dimension.

 \vskip 5mm
In section  4 we focus on the rank $2$ case and we introduce basic
concepts for the rest of the paper. In section 5 we  assume that
the base divisor  $F_V$ is a divisor with normal crossings  and we
find a lower bound for the invariant $\delta (S)$. Section 6 is
devoted to the construction of examples.

\vskip 5mm In sections 7 and 8 we prove Theorem \ref{main}. First
we consider the surfaces with a fibration of higher base genus,
which can be studied directly by analyzing the relative
invariants. Here $K_S^2=16$ and a structure theorem is obtained.
 For the general case we use the geometry of the canonical map, which, due to the previous results,
 factors through the Albanese morphism  composed with the Gauss map.

 \vskip 5mm
 \textbf{Acknowledgements.} We are grateful to J. Amor\'os and F. Campana for
valuable suggestions during the preparation of this paper.

The work was completed during the third author's stay at the
Institut de Matem\`{a}tiques de la Universitat de Barcelona
(IMUB). This stay was supported by the Generalitat de Catalunya's
PIV programm.

\vskip 10mm

\section{Notation}
\vskip 3mm \indent
 All the varieties considered in this paper are reduced, irreducible, projective and defined over the complex numbers.
\vskip 3mm For any coherent sheaf $\mathcal E$ on a variety $X$ we
denote by $\psi_{2,\mathcal E}$ the wedge map
$$\psi_{2,\mathcal E}:\Lambda ^2H^0(X,\mathcal E)\longrightarrow H^0(X,\Lambda ^2\mathcal E).$$
If no confusion arises we will simply write $\psi_2$.

\vskip 3mm Given a coherent sheaf $\mathcal E$ on a variety we
denote by
$$
\delta (\mathcal E)= \frac 12 (c_1(\mathcal E)^2-2c_2(\mathcal
E))\in H^{4}(X,\mathbb Q)
$$ the codimension $2$ piece of the Chern
character. Recall that $\delta $ is additive for exact sequences
and for a locally free sheaf $\mathcal F$, $\delta (\mathcal
F)=\delta (\mathcal F^{\ast})$.

\vskip 3mm We shall say that a cohomological class $\eta \in
H^{2r}(X,\mathbb Q)$ is pseudo effective, and we shall write $\eta
\ge 0$, if
$$
\eta \cdot H_1 \cdots H_{n-r}\ge 0
$$
for any ample divisors $H_1,\dots H_{n-r}.$ The notation
$\eta_1\ge \eta_2$ means $\eta_1-\eta_2\ge 0$.

\vskip 3mm We say that a property holds on $X$ in codimension $k$
if it holds outside of a closed set of codimension at least $k+1$.

\vskip 5mm Throughout the paper $S$ is a smooth minimal surface.

\vskip 3mm We shall say  that $S$  has a fibration of higher base
genus if there exists a fibration $f:S\longrightarrow B$, where
$B$ is a smooth irreducible curve of genus at least $2$.

\vskip 3mm The Albanese map of   $S$  is written as
$a:S\longrightarrow Alb(S).$

\vskip 3mm We denote by $\tau (S)$ the topological index of $S$.
Recall that $\tau(S)=\frac 13(K_S^2-2 c_2(S))$. Hence  $\delta
(\Omega ^1_S)=\frac 32 \tau (S)$.

\vskip 3mm We denote by $F$ the fixed divisor of the canonical
system, and by $M:=K_S-F$ the moving part.

\vskip 3mm The canonical map $S - - \rightarrow \mathbb
{P}(H^0(S,\omega _S)^{\ast })$ induced by $M$ is denoted by
$\varphi_M$ or simply $\varphi $. We put $\Sigma := \varphi_M
(S)$.

\vskip 5mm
\section{On the kernel of $\psi _{2,\mathcal E}$  for a reflexive $\mathcal E$}

The goal of this section is to provide geometrical consequences of
the existence of non trivial elements in the kernel of $\psi
_{2,\mathcal E}$. We are especially interested in the non
decomposable elements. Although the main applications will be for
surfaces and when $\mathcal E$ is the sheaf of differentials, our
aim is to work in a general framework. For this reason we will
start in a quite general setting and we add hypotheses when the
theory needs them to progress.

\vskip 5mm Fix a $n$-dimensional smooth variety $X$ and a
reflexive sheaf $\mathcal E$. In particular, $\mathcal E$ and its
subsheaves are torsion-free. Let $w\in Ker (\psi _{2,\mathcal E})$
and let $V\subset H^0(X,\mathcal E)$ be such that $w\in \Lambda
^2\,V$ and $V$ has minimal dimension with this property. Then $V$
has even dimension $2k$ and we say that this is the rank of $w$.

\vskip 5mm Let us now consider the evaluation map restricted to
$V$:
 $$
 V\otimes \mathcal {O}_X \longrightarrow \mathcal E.
 $$
 The image of this map is a torsion-free sheaf $\mathcal E_V \subset \mathcal
 E$.

 Dualizing the surjective map $V\otimes \mathcal {O}_X\rightarrow \mathcal E_V$,
  we obtain the following short exact sequence
 which serves as definition of the sheaf $N$:
$$
0\longrightarrow \mathcal E^{\ast }_V \longrightarrow V^*\otimes
\mathcal{O}_X \longrightarrow N \longrightarrow 0.
$$

Dualizing again we get an exact sequence
$$
0\longrightarrow N^* \longrightarrow V\otimes \mathcal{O}_X
\stackrel{e}{\longrightarrow} \mathcal E_V^{\ast
\ast}\longrightarrow \mathcal
Ext^1_{\mathcal{O}_X}(N,\mathcal{O}_X) \longrightarrow 0.
$$
Notice that
$$
V\otimes  \mathcal O_X  \longrightarrow  \mathcal E_V^{\ast \ast}
\longrightarrow \mathcal E^{\ast \ast} \cong \mathcal E$$ is again
the evaluation map restricted to $V$, therefore the image of $e$
is $\mathcal E_V\subset \mathcal E_V^{\ast \ast}$. Hence the exact
sequence splits into:
$$
0\longrightarrow N^*\longrightarrow V\otimes \mathcal{O}_X
\longrightarrow \mathcal E_V \longrightarrow 0
$$
and
$$
0\longrightarrow \mathcal E_V\longrightarrow \mathcal E_V^{\ast
\ast} \longrightarrow \mathcal  C_V := \mathcal Ext^1_{\mathcal
O_X}(N,\mathcal O_X) \longrightarrow 0.
$$

\vskip 3mm Observe that contraction with $w$ defines an
isomorphism  $V^{\ast}\otimes \mathcal O_X \longrightarrow
V\otimes \mathcal O_X$. Then the condition $\psi_2(w)=0$
translates into the vanishing of the composition
$$
\mathcal E_V^{\ast }\longrightarrow V^{\ast}\otimes \mathcal O_X
\longrightarrow V\otimes \mathcal O_X \longrightarrow  \mathcal
E_V.
$$
Therefore the following diagram is induced
$$
\begin{CD}
  0 @>>> \mathcal E_V^{\ast }   @>>>  V^*\otimes \mathcal{O}_X @>>> N   @>>> 0 \\
  && @VV\alpha V  @VV\cong V @VV\beta V && \\
  0 @>>> N^{\ast } @>>>  V\otimes \mathcal{O}_X @>>>   \mathcal E_V @>>> 0.
  \end{CD}
$$
Observe (by diagram chasing) that $\alpha $ is injective and
$\beta $ is surjective. Moreover $Coker (\alpha)\cong Ker
(\beta)$. Hence $rank \,\mathcal E _V\le rank \,N$ and
$$rank
\,\mathcal E_V + rank \,N=\dim V=2k.$$ In particular

\begin{lemma} With the above notations we have:
$$
rank \, \mathcal E_V \le k = \frac 12 \dim V.
$$
\end{lemma}

\vskip 5mm The exact sequences above imply several easy
consequences for the invariant $\delta $ of the sheaves involved.
We present them in the following lemma:

\vskip 5mm
\begin{lemma}\label{deltaE}
With the above notations:
\begin{enumerate}
\item [a)]$\delta (\mathcal E_V^{\ast})=\delta (\mathcal E_V^{\ast
\ast}).$
\item [b)]$\delta (\mathcal E_V^{\ast \ast})\ge \delta
(\mathcal E_V).$
 \item [c)]$\delta (\mathcal E_V)+\delta (N^{\ast
})=0.$
\item [d)]$\delta (\mathcal E_V^{\ast})+\delta (N)=0.$
\end{enumerate}
\end{lemma}

\begin{proof}
Since $\mathcal E_V^{\ast}$ is reflexive, then it is locally free
in codimension $2$. Therefore, in the computation of $\delta $ it
can be considered as locally free. This justifies a). Moreover the
sheaf $\mathcal C_V$ is supported in codimension $2$, hence its
restriction to a general surface $S=H_1\cap \dots H_{n-2}$, $H_i$
ample, satisfies $\delta (\mathcal C_{V \mid S})\ge 0$. Therefore
$\delta (\mathcal E_V^{\ast \ast})-\delta (\mathcal
E_V)=\delta(\mathcal C_V)\ge 0$, then b) follows. The rest is
obvious.
\end{proof}

\vskip 5mm
\begin{remark}
We are especially interested in the case where $X=S$ is a surface.
In this situation the inequality of b) becomes an inequality of
integers. Moreover all reflexive sheaves, such as $\mathcal E$,
$\mathcal E_V^{\ast }$ and $N^{\ast }$, are locally free.
\end{remark}

\vskip 5mm \noindent To go further we need the following

\noindent \textbf{Hypothesis:} We assume $rank\, \mathcal E_V=k$.
Notice that this is equivalent to saying that $rank \,N=k$.

With this assumption, the map $$\alpha : \mathcal E_V^{\ast }
\hookrightarrow N^{\ast}$$
  is an injection of reflexive sheaves of the same rank.
  Moreover $c_1(N^{\ast })= -c_1(\mathcal E _V)=c_1(\mathcal E_V^{\ast})$
  ($\mathcal E_V$ is torsion-free, therefore locally free in codimension $1$).
  Hence $\alpha $ is an isomorphism in codimension  $2$. Therefore:
$$
\delta (N^{\ast })=\delta (\mathcal E_V^{\ast })=\delta(\mathcal
E_V^{\ast \ast }).
$$

By implementing  this in \ref{deltaE} we obtain the following
result:

\begin{theorem}\label{thmgen}
Let $X$ be a smooth variety and let $\mathcal E$ be a reflexive
sheaf on $X$. We assume that there exists $w\in Ker
(\psi_{2,\mathcal E})$ of rank $2k$ and let $V\subset
H^0(X,\mathcal E)$ a subspace of dimension $2k$ such that $w\in
\Lambda ^2V$. Let $\mathcal E_V$ be the image of the evaluation
map restricted to $V\otimes \mathcal O_X$.

\noindent If $rank \,\mathcal E_V =k$, then:
\begin{enumerate}
\item [a)] $\delta (\mathcal E_V)+\delta (\mathcal E_V^{\ast
\ast})=0$ \item [b)] $\delta (\mathcal E_V^{\ast \ast}) \ge 0$ and
$\delta (\mathcal {E}_V)\le 0$. \item [c)] Assume that $\mathcal E
_V\hookrightarrow \mathcal E$ is an isomorphism in codimension
$1$. Then $\delta (\mathcal E)\ge 0$
\end{enumerate}
\end{theorem}
\begin{proof}
a) and b) follow from \ref{deltaE}. To see part c) we observe that
now also $\mathcal E_V^{\ast \ast}\hookrightarrow \mathcal E$ is
an isomorphism in codimension $2$, then $\delta (\mathcal E) =
\delta (\mathcal E_V^{\ast \ast})\ge 0$.
\end{proof}

\vskip 5mm In general the class $\delta(\mathcal E)$ is not
positive, even for globally generated sheaves. Consider, for
instance, the symmetric product $S^2(C)$ of a non-hyperelliptic
curve  $C$. Then the Hodge numbers are easily computed and the
equality $\tau(S^2(C))=\delta (\Omega ^1_{S^2(C)})=-g(C)+1$ is
obtained.

Special cases where the inequality $\delta (\mathcal E)\ge 0$ is
true are discussed by Lazarsfeld in \cite{Laz}, ex. 8.3.18.

\vskip 5mm We apply \ref{thmgen}  to generalized Lagrangian
varieties. So we assume $\mathcal E=\Omega ^1_X$. We write
$\overline{\Omega ^1_X}$ instead of $\Omega ^1_{X,V}$. Then we
have:
\begin{corollary}
Let $X$ be a generalized Lagrangian variety such that
$\overline{\Omega ^1_X}\hookrightarrow \Omega ^1_X$ is an
isomorphism in codimension $1$. Then $\delta(\Omega ^1_X)\ge 0$.
\end{corollary}

\vskip 5mm
\begin{remark}\label{Lag}
If $X$ is Lagrangian with $a:X\longrightarrow A$ an immersion,
then we obtain the  well-known isomorphism between the dual of the
tangent and the normal bundle of $X$ in $A$,  providing the
self-dual exact sequence
$$
0\longrightarrow \Omega ^{1 \ast}_X \longrightarrow V^*\otimes
\mathcal{O}_X \longrightarrow \Omega _X^1 \longrightarrow 0.
$$
 Therefore the  pieces of even codimension in the Chern
character of $\Omega ^1_X$ vanish. If $a$ does not have a branch
divisor then $\delta(\Omega ^1_X)\ge 0$ by the last Corollary. In
particular, if $X=S$ is  a Lagrangian surface with these
conditions, then $\tau(S)\ge 0$.
\end{remark}

\vskip 5mm
\section {Rank $2$ vector bundles}

 We apply the results of the preceding section in the
following situation:

We consider $\mathcal E$ a locally free sheaf of rank $2$ with
determinant $L$, and
$$
w=\omega _1 \wedge \omega _2 + \omega _3 \wedge \omega _4
$$
a rank $4$ element in $Ker(\psi_{2,\mathcal E})$. We also assume
that $\omega_1,\dots,\omega_4$ are linearly independent sections
generating a vector space $V\subset H^0(X,\mathcal E)$. The
evaluation map $V\otimes \mathcal O_X\longrightarrow \mathcal E$
is assumed to be generically surjective, in other words $\mathcal
E_V$  has rank $2$.

\vskip 5mm We now introduce some notation. Let us consider the
natural map
$$
\Lambda ^2V\otimes \mathcal O_X \longrightarrow L.
$$
The image  can be written as $L(-F_V)\otimes \mathcal I_{Z_V}$,
where $F_V$ is a divisor and $Z_V$ is a subscheme of codimension
$\ge 2$. In other words, the image of $\Lambda ^2V\longrightarrow
H^0(X,L)$ defines a subsystem of $\vert L\vert$. Then  $F_V$ is
its base divisor and $Z_V$ is the  base locus of its moving part.

\vskip 5mm With this notation, part c) of \ref{thmgen} translates
into
\begin{theorem}\label{thmsurFzero}
With the above notations, assume $F_V=0$. Then:
$$
\delta (\mathcal E) \ge 0.
$$
In particular, if $X=S$ is a generalized Lagrangian surface with
$F_V=0$, then $\delta (\Omega ^1_S)\ge 0$.
\end{theorem}

\vskip 5mm We want to extend this result to Lagrangian surfaces
with some conditions on $F_V$. We will need the following:

\vskip 5mm
\begin{proposition}\label{detomegabar}
The following equality hold:
$$
det\, \overline {\Omega ^1_S}^{\ast \ast }=\omega_S(-F_V).
$$
\end{proposition}

\begin{proof}
 Since
 $$V\otimes \mathcal O_S \twoheadrightarrow \overline {\Omega ^1_S} \subset
{\overline {\Omega ^1_S}}^{\ast \ast }\subset \Omega ^1_S,
$$
and
$$
\overline {\Omega ^1_S} \subset {\overline {\Omega ^1_S}}^{\ast
\ast }
$$
is an isomorphism in codimension $1$, then
$$
\Lambda ^2V \otimes \mathcal O_S \twoheadrightarrow \omega
_S(-F_V)\otimes \mathcal I_{Z_V}\subset \Lambda ^2\overline{\Omega
^1_S}^{\ast \ast }
$$
being the last inclusion an isomorphism in codimension $1$. Hence:
$$
\omega _S(-F_V)=c_1(\omega _S(-F_V)\otimes \mathcal I_{Z_V})=c_1(
\Lambda ^2\overline{\Omega ^1_S}^{\ast \ast })=det\, \overline
{\Omega ^1_S}^{\ast \ast }.
$$
\end{proof}

\vskip 5mm
\section{Bounding $\delta(\Omega ^1_S)$ in the case where $F_V$ is a reduced connected divisor with normal
crossings and rational components}

We want to bound from below $\delta (\Omega ^1_S)$ under some
conditions on $F_V$.  In particular we shall assume that $F_V$ is
a reduced divisor with normal crossings. The principal ingredient
will be a relation of $\overline{\Omega ^1_S}$ with a sheaf of
logarithmic differentials, as we shall see in the proposition
below.

\vskip 3mm We next recall a definition given in the Introduction.

\begin{definition}\label{Fcont}
An effective reduced divisor $D$ on the smooth minimal surface $S$
is contracted by the vector space
 $V\subset H^0(S,\Omega ^1_S)$ if the composition:
$$
V\otimes  \mathcal O_S\longrightarrow \Omega^1_S \longrightarrow
\Omega ^1_{S_{|D}}\longrightarrow \omega
 _D
$$
vanishes.
\end{definition}

\vskip 3mm When $V=H^0(S,\Omega^1_S)$ this composition is the dual
of the differential of the Albanese map $a$ restricted to $D$.
Therefore $D$ is contracted by $H^0(S,\Omega^1_S)$ if and only if
$a(D)$ has dimension $0$.

On the other hand, if the components of $D$ are rational, then, by
the very definition, $D$ is contracted by any $V$.

\vskip 5mm
\begin{proposition}\label{Fsimple}
Assume $F_V$ is a reduced divisor with normal crossings, with
smooth irreducible components which is contracted by $V$. Then
$$
\overline{\Omega ^1_S}^{\ast \ast}\cong \Omega
^1_S(log(F_V))(-F_V).
$$
\end{proposition}

\begin{proof}
Let $p$ be a  smooth point of $F_V$ and let $x,y$ local
coordinates at $p $ such that $y=0$ is the local equation of $F$.
Observe that $\Omega ^1_S(log(F))(-F)$ is generated at $p$ by
$ydx$ and $dy$.

A global $1$-form $\omega \in V$ is locally at $p$: $\omega
=a(x,y)dx+b(x,y)dy$. Since $F_V$ is contracted by $V$, $\omega $
vanishes in the tangent direction $\partial /\partial x $ to
$F_V$, hence $a(x,y)=ya_0(x,y)$.

Then in the complementary of a finite set of points the elements
of $V$ belong to $H^0(S,\Omega ^1_S(log(F_V))(-F_V))$. Since they
are sections of locally free sheaves we conclude that $V\subset
H^0(S,\Omega ^1_S(log(F_V))(-F_V))$.

 The commutativity of the diagram
$$
\begin{CD}
  V\otimes \mathcal O_S & \,\hookrightarrow \, & H^0(S,\Omega ^1_S(log(F_V))(-F_V))\otimes \mathcal O_S \\
  @VVevV  @VVevV  \\
  \overline{\Omega ^1_S} &&    \Omega ^1_S(log(F_V))(-F_V)\\
   \cap && \cap \\
  \overline {\Omega ^1_S}^{\ast \ast}  & \hookrightarrow &  \Omega ^1_S ,
  \end{CD}
$$
and the surjectivity of the left evaluation map gives the
inclusion of subsheaves of $\Omega ^1_S$:
$$
\overline{\Omega ^1_S}^{\ast \ast }\subset \Omega
^1_S(log(F_v))(-F_V).
$$
Since both have the same determinant $\omega _S(-F_V)$ (see
\ref{detomegabar}), they are equal.
\end{proof}

\vskip 5mm
\begin{corollary} \label{norcross}
With the same hypotheses on $F_V$, one has $$\delta (\Omega
^1_S)\ge K_S F_V+\frac 12F_V^2.$$
\end{corollary}
\begin{proof}
We consider the well-known short exact sequence:
$$
0\longrightarrow \Omega ^1_S \longrightarrow \Omega ^1_S
(log(F_V)) \longrightarrow\mathcal O_{F_V} \longrightarrow 0.
$$
Tensoring with $\mathcal O_F(-F)$ and applying \ref{Fsimple} we
obtain:
$$
0\longrightarrow \Omega ^1_S (-F_V)\longrightarrow
\overline{\Omega ^1_S}^{\ast \ast} \longrightarrow\mathcal
O_{F_V}(-F_V) \longrightarrow 0.
$$
Then, by using \ref{thmgen}, b) (recall that here $\mathcal
E_V=\overline{\Omega^1_S}$) and the previous exact sequence, we
have:
$$
\begin{aligned}
0&\le \delta (\overline{\Omega ^1_S}^{\ast \ast })=\delta (\Omega
^1_S(-F_V))+\delta (\mathcal O_{F_V}(-F_V))\,
\\
&=\, \delta (\Omega ^1_S)-K_SF_V+F_V^2+\delta (\mathcal
O_{F_V}(-F_V))=\delta (\Omega ^1_S)-K_SF_V+F_V^2-\frac 32F_V^2\,
\\
& =\,\delta (\Omega ^1_S)-K_SF_V-\frac 12F_V^2.
\end{aligned}
$$
\end{proof}

One of the more interesting consequences of this corollary is the
following result:

\vskip 5mm
\begin{theorem}\label{thmsurFnc}
Assume $F_V$ is a  reduced connected divisor with normal crossings
contracted by $V$. Then $\tau(S)\ge 0$.
\end{theorem}

\begin{proof}
Since $F_V$ is contracted by $V$ we can apply \ref{norcross} to
compute the index:
$$
\tau (S) \ge  \frac 23 K_SF_V+\frac 13 F_V^2.
$$
Moreover $K_S F_V\ge -2-F{_V}^2$. So
$$
\tau (S) \ge  -\frac 43 -\frac 13 F_V^2.
$$
Since $F_V^2\le -2$, then $\tau (S)\ge -\frac 23$. Therefore the
integer $\tau(S)$ is non-negative.
\end{proof}

\vskip 5mm
\begin{remark} \label{remsurFnc}
In fact, it is a straightforward computation to check that the
theorem is true in a more general setting: it is enough to assume
that $F_V$ is reduced contracted by $V$ and it contains at most
one connected component with arithmetic genus $0$. For instance,
if $F$ is a reduced rational connected divisor with normal
crossings, then $F_V$ is contracted and $\tau (S)\ge 0$.

\end{remark}

\vskip 5mm
\section{Examples}

\vskip 5mm The main purpose of this section is to investigate  how
far the topological statement $\tau (S) \ge 0$ remains true when
the condition $F_V=0$ in \ref{thmsurFzero} is weakened. In the
previous section we have seen that $F_V$ rational, connected and
with normal crossings again implies $\tau (S)\ge 0$ (see
\ref{thmsurFnc}).

\vskip 3mm
 Nevertheless, the next example shows a way to construct generalized Lagrangian surfaces
 with negative index.
 In these examples the divisor $F_V$ is not
 contracted by $V$.

\vskip 5mm
\begin{example}
Let  $f:X\longrightarrow Y$ a double covering of a surface $Y$.
Let $C$ be the discriminant divisor.  We assume that $X$, $Y$ and
$C$ are smooth. Then
$$
2 \chi _{top}(Y) - \chi _{top}(C)=\chi _{top}(X).
$$
So $2c_2(Y)-(2-2g(Y))=2c_2(Y)+ C^2+CK_Y=c_2(X).$ On the other hand
$K_X^2=2K_Y^2+2CK_Y+ \frac 12 C^2.$ Together, this gives
$$
\delta (X) = 2 \delta (Y) - \frac 34 C^2.
$$
Now, we choose as $Y$ the product of two curves of genus $\ge 2$,
and $C$ a smooth curve on $Y$ with $C^2>0$. Then $X$ inherits the
two fibrations on $Y$, hence there are four differential forms
$\omega _1,\dots,\omega _4$ such that $\omega _1\wedge \omega
_2=0$ and $\omega _3 \wedge \omega _4=0$. The sum of these two
relations provides a rank $4$ element in $Ker(\psi _{2,\Omega
^1_S})$. But it is easy to check $\delta (Y)=0$ and hence $\delta
(X)<0$.
\end{example}

In the rest of this section we give some examples of generalized
Lagrangian surfaces.

\vskip 5mm \begin{example} The most obvious examples  are surfaces
with two fibrations of higher base genus, for instance those which
are isogenous to a product of curves (see \cite{Cat}). Indeed,
there are two different decomposable elements in the kernel of
$\psi_{2,\Omega^1_S}$ provided by the two fibrations. The addition
gives a rank $4$ element in this kernel.
\end{example}

\vskip 5mm
\begin{example}
Arguing as in Example 6.1 we see that coverings
$f:X\longrightarrow Y$ of generalized Lagrangian surfaces  are
also generalized Lagrangian. Moreover if the covering $f$ is
\'etale then $\tau (X)$ is a multiple of $\tau (Y)$, hence $\tau
(X)\ge 0$ if and only if  $\tau (Y)\ge 0$.
\end{example}

\vskip 5mm
\begin{example}
 In the paper \cite{Bog:Tsc}, Lagrangian surfaces embedded in Abelian four-folds are
 constructed. The authors proved that some of these examples have not fibrations of higher base genus.

 We also recall the important examples  by Campana
  \cite{Cam} and Sommese-Van de Ven \cite{SVdV} of non-trivial
  kernel of the cohomology map. In particular, the paper \cite{Cam}
  revived interest in the nilpotent completion of K\"ahler groups,
  showing examples of exotic nilpotent K\"ahler (conjectured not
  to exist).
  \end{example}

\vskip 5mm
\begin{example}\label{pgmenor4q-11}
If $S$ is a minimal surface with $p_g\le 4q-11$, then the kernel
of $\psi_{2,\Omega ^1_S}$ meets the set of elements in $\Lambda
^2H^0(S,\Omega ^1_S)$ with rank less or equal than $4$. So either
they have a  fibration of higher base genus or they are
generalized Lagrangian. In particular this happens when $p_g=5$
and $q=4$ and we will see that $\tau (S)\ge 0$ in this case.
\end{example}

\vskip 5mm
\begin{example}
We consider the surface constructed in \cite{Pi:Xiao} (pages
85-86) by means of the Prym map of simple cyclic triple coverings
of elliptic curves. Roughly speaking,  the moduli space
$$
  \mathcal {M}:=\{\pi:C\longrightarrow E\,\vert \, g(C)=4, deg\,(\pi )=3,g(E)=1, \pi \text{
   cyclic}\}/\cong
$$
and the map
$$
\begin{aligned}
P: \mathcal M&\longrightarrow \mathcal {A}_3(\gamma  ) \\
\pi &\longmapsto P(C,E)=P(\pi )=J(C)/\pi ^{\ast}(J(E)),
\end{aligned}
$$
are considered. The variety ${A}_3(\gamma )$ is the moduli space
of Abelian threefolds with a polarization of a convenient type
$\gamma $.

Then the general fibre has dimension $1$. In particular there are
$1$-dimensional families of coverings $\{\pi _t\}_{t\in B}$
 with fixed $P(\pi _t)$. Hence there exist a smooth surface $S$,
 a fibration $f:S\longrightarrow B$  and an action of $\mathbb{Z}/(3)$ on
 $S$ which restricts to the generic fibre $C_t$ with quotient an elliptic curve $E_t$.
  Moreover $P(C_t,E_t)$ is a constant polarized Abelian variety
  $P$, which is exactly the kernel of the surjective map $Alb(S)\longrightarrow
  Alb(B)=J(B)$, in particular $q(S)=g(B)+3$.

  By choosing a generator $\rho $ of the group of characters we have a
  decomposition as follows
  $$
  H^0(S,\Omega ^1_S)=\pi ^{\ast }(H^0(B,\omega _B))\oplus
  H^0(S,\Omega ^1_S)^{\rho }\oplus H^0(S,\Omega ^1_S)^{\rho ^2}.
 $$
We can assume, possibly changing the generator, that $dim \,
H^0(S,\Omega ^1_S)^{\rho }=1$ generated by $\alpha $ and $dim \,
H^0(S,\Omega ^1_S)^{\rho ^2}=2$ generated by $\beta _1$ and $\beta
_2$.

Now the equivariant map
$$
\Lambda ^2H^0(S,\Omega ^1_S) \longrightarrow
H^0(S,\omega_S)\longrightarrow H^0(C_t,\omega _{C_t})
$$
sends
$$
\Lambda ^2H^0(S,\Omega ^1_S)^{\rho }=\pi
^{\ast}(H^0(B,\omega_B))\wedge \alpha \oplus \langle \beta _1
\wedge \beta _2\rangle
$$
to the $1$-dimensional space $H^0(C_t,\omega_{C_t})^{\rho}$.
Therefore there exists a $1$-form $\eta $ on $B$ such that the
$2$-form on $S$: $\pi ^{\ast }(\eta )\wedge \alpha + \beta
_1\wedge \beta _2$ vanishes when it restricts to the generic fibre
$C_t$. It is not difficult to see that this produces an element
$w\in Ker (\psi _{2,\Omega ^1_S})$ of rank $4$.

Finally we observe that this element ``does not come from
fibrations", in other words $\beta_1\wedge \beta_2 \ne 0$. Indeed,
otherwise there would exist a fibration of higher base genus
$S\longrightarrow \Gamma $ such that $C_t\longrightarrow \Gamma $
is a covering. Notice that $g(\Gamma )\ne 4$ (otherwise $C_t\cong
\Gamma $ and $f$ would be isotrivial) and $g(\Gamma )\ne 3$ (by
Riemann-Hurwitz). Moreover the existence of $\beta _1,\beta _2$
guarantees that $g(\Gamma )\ge 2$. Hence $g(\Gamma )=2$. This is
impossible since for a generic element $(C_t\longrightarrow E_t)
\in \mathcal M$ the curve $C_t$ does not admit non-trivial maps on
curves of genus $2$.

The previous construction can be used to obtain an actual
Lagrangian surface in $P\times E$, where $E$ is the elliptic
curve, which is the algebraic part of the Intermedian Jacobian
associated to $H^3(P,\mathbb Z)$. We do not know any Lagrangian
subvariety of a simple Abelian variety.
\end{example}

\vskip 5mm
\section{Surfaces of general type with $q=4$, $p_g=5$ and a fibration of higher base genus}

For dimensional reasons, this is the case with the lowest
invariants in which we can ensure the existence of a
non-decomposable 2-form in the kernel of the map

$$ \psi_2:\Lambda ^2 H^0(S, \Omega ^1_S)\longrightarrow H^0(S,\omega_S).
$$

\noindent (at least when there is not a fibration of higher base
genus).

The following is the main result of the rest of the paper

\vskip 5mm
\begin{theorem} \label{indiceq4p5}
Let $S$ be a minimal complex projective surface over $\mathbb{C}$,
of general type, with $q=4$, $p_g=5$. Then $\tau (S) \geq 0$.
\end{theorem}

\vskip 5mm
\begin{remark} Notice that $K_S^2+c_2(S)=24$, therefore $\tau (S)\ge 0$ is equivalent to
$K_S^2\ge 16$. This means that $16 \leq K^2_S \leq 18$, the upper
bound given by the Miyaoka-Yau inequality.
\end{remark}

Observe that, with these invariants, the kernel of $\psi_2$ is
non-trivial. In order to prove the theorem we consider separately
the two natural cases depending on the existence or not of
decomposable elements in the kernel. In this section we give a
complete classification of the first case assuming the existence
of a fibration of higher base genus on $S$, and we postpone the
rest of the proof until the next section.

\vskip 5mm
\begin{theorem}\label{fibred}
The minimal  surfaces of general type with a fibration of higher
base genus and irregularity $4$ and geometric genus $5$ are
isotrivial, of the form $(C\times H)/\mathbb Z_2$, where:
\begin{itemize}
\item[a)] $C$ and $H$  are two curves of genus $3$ equipped with
involutions without fixed points, or

\item[b)]$C$ is a curve of genus $5$ equipped with an involution
without fixed points and $H$ is a bielliptic curve of genus $2$.
\end{itemize}
In both cases $\mathbb Z_2$ acts diagonally.

Both examples give irreducible families of surfaces of dimensions
$6$ and  $8$ respectively.

All of them satisfy $K_S^2= 16$. In particular $\tau (S)=0$.
\end{theorem}

\begin{proof}

Let $f:S\longrightarrow B$ be a fibration with $g(B)=b\ge 2$ and
fibre $H$ of genus $g\ge 2$. We consider the relative invariant
$$0\le \chi _f =\chi (\mathcal{O}_S)-(b-1)(g-1)=2-(b-1)(g-1).$$
There are two cases: \vskip 3mm 1. $b=g=2$, then $q(S)=4=b+g$ and
so $S\cong B\times H$ which contradicts $p_g(S)=5$. \vskip 3mm 2.
$b=2,g=3$ or $b=3,g=2$. Therefore $\chi_f=0$ and then also
$0=K_{S\vert B}^2=K_S^2-8(g-1)(b-1)$. So $K_S^2=16$.

\vskip 3mm We can go further in the analysis of the fibrations in
case 2. First we recall that
 $\chi_f=0$ implies $f$ isotrivial and smooth. Then we can apply the
  general results on these surfaces and, in particular, the very explicit
  description provided in \cite{Se}.  For instance we know that $S$ is the quotient
  of a product of two smooth curves $C\times H$ by a finite  group $G$ of order $n>1$
  acting diagonally and the map $C\longrightarrow C/G$ is unramified
  (since $f$ is smooth). Then there  are two fibrations $g_i: S\longrightarrow D_i/G$
  which are the composition of
$$S\cong (C\times H) /G \longrightarrow C/G\times H/G $$ with the
projections. The initial fibration $f$ equals $g_1$, so $B=C/G$.

\vskip 3mm Observe that $$
\begin{aligned}
q(S)\,=&\,q((C\times H)/G)=h^0((C\times H)/G,\Omega ^1_{(C\times
H)/G})\,=
\\
&\dim H^0(C\times H,\Omega ^1_{C\times
H})^G=g(C/G)+g(H/G)=b+g(H/G),
\end{aligned}
$$ hence $4=b+g(H/G)$.

\vskip 5mm Assume first $b=g(H/G)=2$. Then, Riemann-Hurwitz for
$H\longrightarrow H/G$ implies $n=2$. On the other hand $$
\begin{aligned}
&p_g(S)=p_g((C\times H)/G)=h^0((C\times H)/G,\omega _{(C\times
H)/G})\,=
\\
&\dim H^0(C,\omega _{C})^+ \cdot \dim H^0(H,\omega _{H})^++ \dim
H^0(C,\omega _{C})^-\cdot \dim H^0(H,\omega _{H})^-
\\
&=\,b \cdot g(H/G)+(g(C)-b)(g(H)-g(H/G))=4+(g(C)-2).
\end{aligned}
$$ Therefore $g(C)=3$ and $C\longrightarrow B$ is unramified.

\vskip 5mm Assume now $b=3$ (hence $g=2$), $g(H/G)=1$. In this
situation $S\longrightarrow H/G$ has two singular fibres since
$H\longrightarrow H/G$ has two ramification points. These fibres
are of the form $(\# H_i)C/H_i$, $i=1,2$, where $H_i$ are
subgroups of $G$. Then, by using the formula of  theorem 4.1 in
\cite{Se} for $K_S$, $$
\begin{aligned}
&2g(H)-2=2=H^2+K_SH=K_SH=
\\
&(\#H_1-1)H\cdot (C/H_1)+ (\#H_2-1)H\cdot (C/H_2).
\end{aligned}
$$ Since $H\cdot (C/H_i)\ne 0$, we deduce $H\cdot (C/H_i)=1$ and $\#
H_i=2$. On the other hand $$ n=HC=(\# H_1)H\cdot C/H_1=2. $$ Hence
$G\cong \mathbb Z_2$. Computing $p_g(S)$ as above we find $g(C)=5$
and that $C\longrightarrow B$ is unramified.
\end{proof}

\vskip 5mm Notice that the previous considerations allow us to
prove that the  Albanese dimension  is always $2$. Indeed, $S$
would have otherwise  a fibration on a curve of genus 4, which is
not compatible with the range of genus found above.

\vskip 5mm \section{Surfaces of general type with $q=4$, $p_g=5$
and without a fibration of higher base genus}

\vskip 5mm
 In this section we finish the proof of \ref{indiceq4p5}. We
 assume now that $S$
 has not a fibration of higher base genus and the dimension of $a(S)$ is $2$. In particular,
 they are generalized Lagrangian (see \ref{pgmenor4q-11}). There exists a
non-decomposable $$w\in \Lambda ^2 H^0(S,\Omega ^1_ S),$$ unique,
up to constant, which vanishes in $H^0(S,\omega _S)$. Then we can
use the notations and theorems of \S 4. Here $V=H^0(S,\Omega
^1_S)$, $F_V=F$ the base divisor of the canonical system and $\psi
_2$ is surjective. Also by \ref{thmsurFzero} we can restrict
ourselves to the hypothesis $F\ne 0$. In particular, by
$2$-connectivity, $MF\ge 2$, where $M$ is the moving part of the
canonical system.

\vskip 5mm

Now we analyze the canonical map of $S$. Note that a general type
surface with $q \geq 3$ has a canonical image of dimension 2
(\cite{Xiao2}).

Observe we have a short exact sequence of vector spaces $$
0\longrightarrow \left\langle w \right\rangle \longrightarrow
\Lambda ^2  V \longrightarrow H^0(S,\omega_S) \longrightarrow 0.
$$ By dualizing, we see $\mathbb P(H^0(S,\omega _S)^{\ast })$ as
the hyperplane $H_w\subset \mathbb P(\Lambda ^2 V^{\ast })$
naturally attached to $w$. Notice that $H_w$ is not tangent to the
Grassmannian $Gr(2,V^{\ast})$, since $w$ is not decomposable.

\vskip 5mm  We attach to a point in $S$ the geometric fibre of
$\overline {\Omega ^1_S}\,^{\ast}\subset V^{\ast }\otimes \mathcal
O_S$ at this point. We see this $2$-dimensional vector space as an
element in the Grassmannian naturally embedded in $\mathbb
P(\Lambda ^2V^{\ast })$. By construction we have a commutative
diagram of rational maps: $$
\begin{CD}
  S @>>> Gr(2,V^{\ast} ) \\
  @V\varphi VV  @VVV  \\
  H_w @>>>  \mathbb P(\Lambda ^2 V^{\ast }).
  \end{CD}
$$
As a consequence we obtain that the image $\Sigma $ of the
canonical map is contained in the smooth quadric
$Q:=Gr(2,V^{\ast})\cap H_W$ of $\mathbb P^4=H_W=\mathbb
P(H^0(S,\omega _S)^{\ast})$.

Any divisor on a smooth quadric of $\mathbb P^4$ is a complete
intersection, so there exists a hypersurface $T$ of degree $t$
such that $\Sigma = Q\cap T$. Denote $e=\text {deg}\,(\varphi )$.
We have $$ M ^2=\text{deg}\,(\Sigma )\,\text {deg}\,(\varphi
)+b=2\,t\,e\,+\,b.
$$ where $b$ is the contribution of the base points of the linear system $|M|$. The surface $\Sigma $ is not
degenerated by construction, hence $t$ is at least $2$.

\vskip 5mm
\begin{proposition} \label{casosbajos}
Under our hypothesis
\begin{itemize}
\item[a)] $M^2 \geq 10$. \item[b)] If $M^2=10,11$, then $b=0,1$
respectively.
\end{itemize}
\end{proposition}

\begin{proof}
If $te \geq 5$ we are done. We assume $t\, e \,\le 4.$ There are
the following possibilities:

\vskip 3mm Case A: \quad $t= 2$,  $e=2$.

Here $S$ is a double covering of the intersection of two quadrics
$Q, T$ in $\mathbb P^4$. Since one of them is smooth it is a
(possibly singular) del Pezzo surface. Hence $S$ is covered by a
linear system of hyperelliptic curves. This contradicts the
results of \cite{Xiao} (see also \cite{Pi:Du}).

\vskip 3mm Case B: \quad $t= 3$,  $e=1$.

In this case the canonical map is a birational morphism with a
surface $\Sigma $ which is not of general type; this contradicts
our hypothesis.

\vskip 3mm Case C: \quad $t= 4$,  $e=1$.

We consider a general hyperplane section of $\Sigma $. It is the
complete intersection of a smooth quadric and a quartic surface in
$\mathbb P^3$. So its arithmetical genus is $9$. Its preimage is a
general curve $C$ in $|M|$. We have $M^2+MK_S\ge 16+MF\ge 18$. If
$C$ is smooth, then genus $C$ is at least $10$, and it is a
desingularization of its image, $\varphi$ being birational, a
contradiction. If the general element $C$ is singular, then by
Bertini's theorem it must be singular at the base points of the
linear system. In this case, necessarily $b \geq 4$ and so
$M^2\geq 12$.
\end{proof}

Now we study some properties of special curves on $S$.

\vskip 5mm
\begin{proposition} \label{curvasespeciales}
Under our hypotheses
\begin{enumerate}
\item[a)] There is not an effective divisor $D$ on $S$ with
$D^2=0$ and $K_SD\leq 2$. \item[b)] If $C$ is a rational curve on
$S$, then $C\subseteq F$. \item[c)] Let $k$ be the number of
(-2)-curves on $S$. Then $$K^2_S \leq 18-\frac{3}{4}k.$$ If
equality holds then the (-2)-curves are disjoint.
\end{enumerate}
\end{proposition}

\begin{proof}

a) The arithmetical genus of $D$ is at most $2$. Let $\Gamma
=a(D)$ be the image of $D$ by the Albanese map. Since the square
of $D$ is not negative it cannot be contracted, hence $\Gamma $ is
a curve with arithmetical genus $1$ or $2$.  By definition the
image of $D$ by the map $a':S\longrightarrow Alb(S)/\left\langle
\Gamma \right\rangle $ is a point, hence, since $D^2=0$, $a'(S)$
is a curve, and generates the Abelian variety in the image, so its
genus is at least 2. Then $a'$ defines a fibration of higher base
genus and this contradicts our hypothesis on $S$.

\vskip 3mm

\noindent b) The global differential $1$-forms vanish along the
tangent directions to a rational curve. Then the kernel of the
form at each point of the curve is $1$-dimensional and the wedge
product of two of them is zero. In other words, a rational curve
is contained in the base locus of $\Lambda ^2H^0(S,\Omega
^1_S)=H^0(S,\omega _S)$, which is $F$.

\vskip 3mm

\noindent c) The result follows from \cite{Miy}, Theorem 1.1,
taking $D=0$ and E the divisor of (-2) curves on $S$. Let
$E=E_1+...+E_l$ be its decomposition in connected components. Let
$k_i$ be the number of irreducible components of $E_i$ (we have
then that $k=k_1+...k_l$). Then the Zariski decomposition of
$K_S+E$ is exactly $P=K_S$, $N=E$. For any $i$ define

$$\nu _i=e(E_i)-\frac{1}{|G_i|}$$

where $G_i$ is the finite group of the quotient singularity
produced contracting $E_i$. Clearly $e(E_i)=1+k_i$ and $|G_i|\geq
2$ (and when equality holds $E_i$ is a single (-2)-curve). Then,
applying Theorem 1.1 in \cite{Miy}, we obtain

$$\frac{1}{2}l+k \leq \Sigma_i \nu _i \leq
c_2(S)-\frac{1}{3}K_S^2-\frac{1}{4}E^2=24-\frac{4}{3}K_S^2+\frac{1}{2}l$$

and so

$$K_S^2 \leq 18 - \frac{3}{4}k.$$

If equality holds, then $|G_i|=2$ for all $i$ and so all the
connected components of $E$ are single (-2)-curves.

\end{proof}

In the next result we present useful information on the base locus
of the linear system $\vert K_S+M \vert$.

\vskip 5mm
\begin{proposition} \label{basepoints}
Under our hypotheses, assume $K_S^2 \leq 15$. Then
\begin{itemize}

\item[a)] If $K_SF=1$, $F^2=-1$, then $|K_S \, + \, M |$ has no
base point. \item[b)] Let $p$ be a base point of $|K_S \, + \, M
|$. Then, there exists an effective divisor $C$ through $p$, with
$C^2=-1 \, , MC=0 \, , FC=1 $. If $KF=0$ we can assume $C$ is an
elliptic smooth curve. If, moreover, $F^2=-2$, then such a curve
$C$ is unique with the above numerical properties. \item[c)] The
base divisor $B$ of $|K_S \, + \, M |$ is formed by elliptic
curves in $F$ contracted by $|M|$. In particular, since $MF\ge 2$,
$F\ne B$.

\end{itemize}
\end{proposition}

\begin{proof}
From the previous results, we know that $M^2\geq 10$, $MF \geq 2$.

Let $p$ be a base point of $|K_S+M|$. Since $M$ is big and nef we
can apply Reider's theorem and obtain an effective divisor $C$
passing through $p$ with one of the following numerical
properties:

\begin{itemize}
\item $C^2=0$, $MC=1$ \item $C^2=-1$, $MC=0$
\end{itemize}

Since M is big, the Hodge-Index theorem says

$$det \begin{pmatrix} M^2 & MF & MC \\ MF & F^2 & x \\ MC & x & C^2
\end{pmatrix} \geq 0 $$

\noindent where $x=FC$. Since $1+x=KC+C^2$ is even, $x$ must be
odd and, since $KC=MC+x \geq 0$ we obtain that $x \geq -MC$.
Observe that, in fact, $x \geq 1$. Indeed, if $x=-1$, then $KC=0$
and so $C$ is a chain of (-2)-curves, which is impossible since
$C^2=0$.

Moreover we have

$$12 \leq M^2 + MF \leq K_S^2 \leq 15$$

\noindent and so $MF \leq 5$. Computing the determinant, we obtain
in the first case

$$-x^2M^2+2xMF-F^2\geq 0$$

\noindent which gives

$$10x^2+2 \leq (x^2-1)M^2+12 \leq 2(x+1)MF \leq 10(x+1)$$

As $x$ is odd and positive, this can only happen if $x=1$. But
then observe that $KC=2$ and $C^2=0$, which contradicts
proposition \ref{curvasespeciales} a).

So, only the second numerical possibility for $C$ may occur. In
this case, the determinant above gives

$$-F^2M^2-M^2x^2+(MF)^2\geq 0$$

\noindent and so,

$$F^2+x^2 \leq \frac{(MF)^2}{M^2}$$

\noindent Hence

$$-MF+x^2 \leq F^2+x^2 \leq 2$$

\noindent which produces $x^2\leq 2+MF \leq 7$ and so $x=1$.

If $KF=1$, $F^2=-1$ we have $(F+C)^2=0$ and $K(F+C)=2$ which again
contradicts proposition \ref{curvasespeciales}. This proves a).

If $KF=0$, let $C_0$ be an irreducible component of $C$ meeting
$F$ ($x=FC=1$). Since $M$ is nef, we still have $MC_0=0$. Since
$C_0$ is irreducible and $M$ big we obtain $-2 \leq C_0^2 \leq -1$
by the Hodge-Index theorem. If $C_0^2=-2$ then $C_0$ would be
rational and so, by \ref{curvasespeciales}, we would obtain $C_0
\leq F$. But this is impossible, since $1=KC_0\leq KF=0$. So $C_0$
is a smooth elliptic curve with $MC_0=0$ and $C_0^2=-1$. Finally,
note that  $C_0$ also contains necessarily $p$.

If, moreover, $F^2=-2$ and we had two curves $C_i$ with
$C_i^2=-1$, $MC_i=0$ and $FC_i=1$ note that, being orthogonal to
$M$, the Hodge-Index theorem would give

$$det \begin{pmatrix} C_1^2 & C_1C_2 \\ C_1 C_2 & C_2^2
\end{pmatrix}> 0$$

\noindent and hence $C_1C_2=0$. So

$$(F+C_1+C_2)^2=0$$ $$K(F+C_1+C_2)=2$$

\noindent which contradicts proposition \ref{curvasespeciales}.
This finishes the proof of b).

In order to prove c) observe that, since $|M|$ has no base
divisor, the base divisor $B$ of $|K_S+M|$ must be contained in
$F$. Let $D$ be an irreducible component of $B$. For any point $p
\in D$ there must exist an effective divisor $C_p$ containing it,
such that $MC_p=0$ by b). We can assume $C_p$ is irreducible with
this property. Such a $C_p$ cannot move with $p$ otherwise $S$
would be covered by contracted curves. So $D=C_p$ for any $p$. But
then

$$0\leq MD\leq MC_p=0$$

\noindent and so $MB=0$. Since $MF \geq 2$, we obtain that
necessarily $B \neq F$ which proves c).
\end{proof}

\vskip 10mm We can finally prove the main theorem of this part:

\vskip 5mm
\begin{proof}(of Theorem \ref{indiceq4p5})

If $S$ has a fibration of higher base genus, then the result is
given in theorem \ref{fibred}. Otherwise, by the previous
considerations, there exists an element $ \omega _1\wedge
\omega_2+\omega _3\wedge \omega_4=0$.

Since $K_S$ is nef and by $2$-connectivity we obtain
$$ K_S^2=K_S (M+F)\,\ge \,K_S M=M^2+MF \,\ge \, M^2+2, $$
therefore $M^2  \,\ge \, 14 $ implies the theorem.

Moreover proposition \ref{casosbajos} gives $M^2 \geq 10$. We have
the following numerical facts:

\begin{enumerate}
\item $M  F$ is even. \item $\chi (S, \mathcal
O_S(M+K_S))=h^0(S,\mathcal O_S(M+K_S))=2+M^2+\frac 12 M F$. \item
$\chi (S, \mathcal O_S(2K_S))=h^0(S,\mathcal O_S(2K_S))=2+K_S^2$
\item $h^0(S, \mathcal O_S(2M))\ge 13$. If $t\ge 3$, then $h^0(S,
\mathcal O_S(2M))\ge 14$ \item $h^0(S,2M)< h^0(S, K_S+M)$
\end{enumerate}

Indeed, (1) follows from adjunction formula for $F$, and (2) and
(3) from Riemann-Roch formula and Ramanujam vanishing theorem. (5)
is exactly \ref{basepoints} (iii). (4) is a consequence of
counting quadrics containing $\Sigma$: there are 2 linearly
independent quadrics through $\Sigma $ if $t=2$, and only 1
otherwise.

All these numerical restrictions configure the following
possibilities if $K^2_S \leq 15$:

\vskip 5mm
\begin{tabular}{|c|c|c|c|c|c|c|c|}
\hline $M^2$ & $MF$ & $K_S^2$ & $F^2$ & $K_SF$ & $h^0(2M)$ &
$h^0(M+K_S)$ & $h^0(2K_S)$   \\ \hline 11 & 4 & 15 & -4 & 0 & 14 &
15
& 17   \\ \hline 12 & 2 & 14 & -2 & 0 & $\geq$ 13 & 15 & 16   \\
\hline 12 & 2 & 15 & -1 & 1 & $\ge $13 & 15 & 17   \\ \hline 13 &
2 & 15 & -2 & 0 & $\ge $13 &16 & 17   \\ \hline
\end{tabular}

\vskip 5mm

Notice that, following proposition \ref{basepoints}, in all these
cases either there are no base points of $|K_S+M|$ (when $KF=1$)
or there are no base divisor (when $KF=0$ all the components of
$F$ are (-2)curves).

We will prove, by a case by case consideration, that none of these
four possibilities may occur.

\vskip 5mm

\noindent Case 1. $M^2=11$, $K^2_S \leq 15$ is not possible.

\vskip 3mm In this case $b=1$ and $te=5$. Let $p$ be the unique
base point of $|M|$. It is necessarily a simple base point. On the
other hand, $t=5$ implies that $\Sigma$ is contained in a unique
quadric of $\mathbb{P}^4$ and so the map
$$S^2H^0(S,M) \longrightarrow H^0(S,2M)$$
is surjective, by dimension counting. Hence $p$ is a singular base
point of the linear system $|2M|$.

Consider the inclusions
$$H^0(S,2M)\subseteq H^0(S,{\mathcal I}_p(K_S+M)) \subseteq
H^0(S,K_S+M).$$

Note that the first subspace has dimension 14 and the last 15. So
the subspace in the middle coincides with one of them. If
$H^0(S,\mathcal{I}_p(K_S+M))= H^0(S,K_S+M)$, then $p$ is a base
point of $|K_S+M|$. By proposition \ref{basepoints} there exists
an elliptic smooth curve $C$ through $p$ with $MC=0$. Since $p$ is
a base point of $|M|$ and it has no base component, this is not
possible.

\vskip 5mm

If $H^0(S,2M)=H^0(S,{\mathcal I}_p(K_S+M))$, then $p$ is not a
base point of $|K_S+M|$ and all the divisors of this linear system
passing through $p$ are singular at $p$. We can apply the second
part of Reider's theorem: given any vector $v \in T_pS$, there
exists a divisor $D$ through $p$, such that $v \in T_pD$ and $MD
\leq 2$.

Assume first that $D$ is smooth at $p$ and, in particular, there
is only one irreducible component $D_v$ of $D$ passing through
$p$. Since it must have $v$ as a tangent vector, and it is
non-singular at $p$, this component moves with $v$. But this is
impossible, since $MD_v \leq 2$ and so, $p$ being a simple base
point of $|M|$, either the curves $D_v$ are contracted via the
canonical map or they are rational curves.

Therefore we can assume that $D$ is singular at $p$. In particular
$MD\ge 2$ and hence $MD=2$. Looking at the possibilities of
Reider's theorem we have $D^2=0$.

Now we compute $x=FD$ by using the Hodge-Index theorem. We
consider the determinant of the matrix given by the intersection
products of $M,F$ and $D$:
$$
det\, \begin{pmatrix}
  11 & 4 & 1 \\
  4 & -4 & x \\
  1 & x & 0 \\
\end{pmatrix}
\ge 0
$$
and we obtain $x=0,1,2$. Notice that $KD+D^2=MD+FD=2+x$ is even,
so $x$ is either $0$ or $2$. Therefore $p_a(D)$ is either $2$ or
$3$ and considering the quotient  $Alb(S)/< a(D)> $, we construct
a fibration $f:S\longrightarrow B$, $g(B)=1,2$ with fibre a
multiple of $D$. Since $F$ consists of rational curves, then
$f(D)=0$. Hence $x=FD=0$. In particular $KD=2$ (and $D^2=0$). By
\ref{curvasespeciales}, part a) we obtain a contradiction.

\vskip 5mm \noindent Case 2. $M^2=12$, $K^2_S \leq 15$ is not
possible.

\vskip 3mm

Assume first $|K_S+M|$ has base points. Then, by proposition
\ref{basepoints} (i) we necessarily have $MF=2 \, , K^2_S=14 \, ,
KF=0$.

Observe that $|K_S+M|$ has at least one base point on $F$. Indeed,
if $|M|$ is base point free so it is $|2M|$ and then this is
immediate. If $b\neq 0$, then $t\neq 2$ and hence $h^0(S,2M)=14$.
Since $h^0(S,K_S+M)=15$ we have that $H^0(S,2M))+F\subseteq
H^0(S,K_S+M)$ is a hyperplane. From $(K_S+M)F=2$ we obtain some
base point $p$ on $F$.

By proposition \ref{basepoints} there exists an elliptic curve
through $p$ with $FC=1$. Moreover, this $C$ is unique so there are
no more base points on $F$. Let $F_0$ be the component of $F$
containing $p$.

Notice that $MF_0=0$ implies that $F_0+C$ satisfies the same
numerical conditions as $C$ and this contradicts proposition
\ref{basepoints} (ii). If $MF_0=2$, then $$
0=K_SF_0=MF_0+FF_0=2+FF_0.$$ The classification of \cite{BPV}
gives $F=F_0$ is rational and irreducible, in particular we can
use remark \ref{remsurFnc} and we obtain $\tau (S)\ge 0$.

Finally assume $MF_0=1$. Observe that $F_0$ has to be a simple
component, since $CF_0=CF=1$. The equality $MF=2$ implies the
existence of another component $F_1$ with $MF_1$. Then
$h^0(S,M+K_S-F_1)=14$ ($\vert M+K_S\vert$ does not have base
divisor) and therefore there is a new base point on $F_1$.
Applying again proposition \ref{basepoints} we find an elliptic
curve $C'$ through this point. Since $C\ne C'$ this contradicts
that $C$ is unique.

Then $FF_0=-1$. We use again the classification in \cite{BPV} to
obtain that $F$ is reduced and with normal crossings. By using
\ref{remsurFnc} we obtain  $\tau (S)\ge 0$.

\vskip 5mm
 Assume now that $\vert M+K_S\vert$ is base point free.
Since $(M+K_S)F\neq 0$, then $\varphi_{M+K_S}(F)$ is not a point.
The dimension of the linear subspace generated by
$\varphi_{M+K_S}(F)$ is precisely $h^0(S,K_S+M)-h^0(S,2M)-1$.
Hence, necessarily $h^0(S,\mathcal O_S(2M))=13$ and therefore
$\varphi_{M+K_S}(F)$ is a line $l$.

\vskip 5mm We decompose $F=F_0+F'$, $F_0$ consisting in the
components of $F$ not contracted by $M$; observe that $F_0$ has at
most two components. Then, the restriction to $F_0$ of the map
$\varphi_{M+K_S}$ gives $h:F_0\longrightarrow l$ of degree
$(M+K_S)F_0\ge 2$. We now apply the second part of Reider's
theorem to the general fibre of $h$: there exists an effective
divisor $E$ such that either $ME=0,\, E^2=-1,-2$ or $ME=1,\,-1,0$
or $ME=2,\,E^2=0$.

Assume that $ME$ is either $0$ or $1$. Then if $E$ moves with the
fibres of $h$ the surface is covered by rational or contracted
curves: a contradiction. Therefore $E\ge F_0$ (if $F_0$ has two
components we apply Reider's theorem to couples of points in
different components), hence $ME\ge MF_0=2$. We deduce that only
the case $ME=2,\,E^2=0$ is possible.

Observe that $EF=0$. Indeed, denote $x=EF$. Since
$E^2+K_SE=ME+x=2+x$ is even, we obtain that $x$ is even. Moreover
the determinant of the matrix of products of $M,F,E$: $$
\begin{pmatrix}
12 & 2 & 2 \\ 2 & F^2 & x \\ 2 & x & 0
\end{pmatrix}
$$ is positive, hence $8x-4F^2-12x^2>0$. Therefore $x=0$.

Now $K_S E=ME=2$ and $E^2=0$ implies the existence of a fibration
of higher base genus (see proposition \ref{curvasespeciales}), a
contradiction.

\vskip 5mm

\noindent Case 3. $M^2=13$, $K^2_S \leq 15$ is not possible.

\vskip 3mm

Since $KF=0$, $F_{red}$ is exactly the divisor of (-2)-curves on
$S$. If $k$ is the number of such curves, and we have $K^2_S=15$,
proposition \ref{curvasespeciales} c) says that $k \leq 3$ (if
$k=4$, then we have equality and then the 4 (-2)-curves are
disjoint: this is impossible if $F^2=-2$). But then $F^2=-2$
implies that $F=A_k$ and corollary \ref{thmsurFnc} finishes the
proof.
\end{proof}


\begin{thebibliography}{999999}


\bibitem{Am:wb}
Amor\'os, J.; Burger, M.; Corlette, K.; Kotschick,D.;Toledo,D.
 {\it Fundamental Groups of Compact K\"ahler Manifolds}.   Mathematical Surveys and Monographs,
 {\bf 44},1996.

\bibitem{BPV}
Barth, W.; Peters, C.;  Van de Ven,A.
    {\it Compact Complex Surfaces}.
   Ergebnisse der Mathematik (3), Springer-Verlag, Berlin,
 1984.

\bibitem{Cat}
    Catanese, F.
    {\it Fibred surfaces, varieties isogenous to a product and related moduli spaces}.
   American J. of Math.  {\bf 122} (2000),1--44.

\bibitem{Cat2}
    Catanese, F.
    {\it Moduli and classification of irregular {K}aehler manifolds (and algebraic varieties) with {A}lbanese general
    type fibrations},
   Invent. Math. {\bf 104} (1991),263--289.

\bibitem{Cam}
    Campana, F.
    {\it Remarques sur les groupes de {K}\"ahler nilpotents}.
   Ann. Sci. ENS {\bf 28} (1995), 307--316.

\bibitem{Miy}
     Miyaoka,Y.
   {\it The maximal number of quotient singularities on surfaces with given numerical invariants}.
  Math. Ann. {\bf 268} (1984), 159--171.

\bibitem{Morgan}
     Morgan,J.
   {\it The Algebraic Topology of smooth algebraic varieties}.
  Publ. IHES {\bf 48} (1978), 137--204.

\bibitem{Se}
     Serrano,F.
   {\it Isotrivial fibred surfaces}.
  Ann. Mat. Pura Appl. {\bf 171} (1996), 63--81.

\bibitem{SVdV}
    Sommese, A.J.; Van de Ven, A.
    {\it Homotopy groups of pullbacks of varieties}.
   Nagoya Math. J. {\bf 102} (1986), 79--90.

\bibitem{Xiao}
     Xiao,G.
   {\it Irregular families of hyperelliptic curves}.
  Algebraic Geometry and Algebraic Number Theory (Tianjin)
   (1989-1990),
  152--156.


\bibitem{Xiao2}
     Xiao,G.
   {\it L'irr\'egularit\'e des surfaces de type g\'en\'eral dont le syst\`eme canonique est compos\'e d'un pinceau}.
  Compositio Math. {\bf 56}, (1985), 251--257.


\bibitem{Reider}
     Reider,I.
   {\it Vector  bundles of rank 2 and linear systems on algebraic surfaces}.
  Annals of Mathematics {\bf 127} (1988),309--316.


\bibitem{Pi:Du}
     Pirola,G.P.
   {\it Curves on generic {K}ummer varieties}.
   Duke Math. J. {\bf 59}
       (1989), 701--708.


\bibitem{Pi:Xiao}
     Pirola,G.P.
   {\it On a conjecture of {X}iao}.
  J. Reine Angew. Math.
       {\bf 431}
     (1992), 75--89.


\bibitem{Bog:Tsc}
    Bogomolov,F.; Tschinkel,Y.
    {\it Lagrangian subvarieties of {A}belian fourfolds}.
   Asian J. Math. {\bf 4} (2000), 19--36.

\bibitem{Laz}
    Lazarsfeld,R.
    {\it Positivity in {A}lgebraic {G}eometry, {I}{I}}.
   Ergebnisse Der Mathematik (3), Springer-Verlag, Berlin,
    2005.


\end{thebibliography}

\bibliographystyle{amsalpha}

Miguel Angel Barja

Departament de Matem\`atica Aplicada I

Universitat Polit\`ecnica de Catalunya

ETSEIB Avda. Diagonal 647

08028 Barcelona (Spain)

Miguel.Angel.Barja@upc.edu \vskip 3mm

Juan Carlos Naranjo

Departament d'\`Algebra i Geometria

Facultat de Matem…tiques

Universitat de Barcelona

Gran Via, 585

08007 Barcelona (Spain)

jcnaranjo@ub.edu \vskip 3mm

Gian Pietro Pirola

Dipartimento di Matematica

Universit\`a di Pavia

Via Ferrata, 1

27100 Pavia (Italy)

pirola@dimat.unipv.it

\end{document}